\magnification\magstep1
\baselineskip = 18pt
\def \Bbb {\bf}
\def\square{\vcenter{\hrule height1pt
\hbox{\vrule width1pt height4pt \kern4pt
\vrule width1pt}
\hrule height1pt}}

\centerline{\bf Volume ratios and a reverse isoperimetric
inequality}\bigskip
 \centerline{Keith Ball$^{(1)}$}
\centerline{Trinity College}
\centerline{Cambridge}\smallskip
\centerline{and}\smallskip
\centerline{Texas A\&M University}
\centerline{College Station, Texas}\bigskip

\noindent {\bf Abstract.} It is shown that if $C$ is an $n$-dimensional
convex body then there is an affine image $\widetilde C$ of $C$ for which

$${|\partial \widetilde C|\over |\widetilde C|^{n-1\over n}}$$

\noindent is no larger than the corresponding expression for a regular
$n$-dimensional ``tetrahedron''. It is also shown that among
$n$-dimensional subspaces of $L_p$ (for each $p\in [1,\infty]),
\ell^n_p$ has maximal volume ratio.\vskip3in

\noindent A.M.S. (1980) Subject Classification: \ 52A20

\noindent $^{(1)}$Supported in part by N.S.F. DMS-8807243
\vfill\eject

\noindent {\bf \S 0. Introduction.}

The famous isoperimetric inequality in ${\Bbb R}^n$ states that among
bodies of a given volume, the Euclidean  balls have least surface area.
Measurable sets of finite volume may have infinite ``surface area'' and, if
$n\ge 2$, even convex bodies of a given volume may have arbitrarily large
surface area if they are very flat. Nevertheless, classical inequalities
such as the isoperimetric inequality do admit of reverse forms: \ the
important reverse Santalo and reverse Brunn-Minkowski inequalities of
[B-M] and [M] are examples.

Probably the most natural way to reverse the isoperimetric inequality is to
consider classes of affinely equivalent convex bodies rather than
individual bodies. The inequality between volume and surface area is proved
only for one representative of each class; the ``least flat'' member of
that class.

Modulo affine transformations it will be shown that among all convex bodies
in ${\Bbb R}^n$, the $n$-dimensional tetrahedron has ``largest'' surface
area for a given volume, while among symmetric convex bodies, the cube is
extremal. More precisely, the principal theorems proved in this paper are
the following: \ (volume and area are denoted throughout by $|\cdot |$).
\medskip

\noindent {\bf Theorem 1.} Let $C$ be a convex body in ${\Bbb R}^n$ and $T$
a regular $n$-dimensional tetrahedron (solid simplex). Then there is an
affine image $\widetilde C$ of $C$ satisfying

$$\eqalign{|\widetilde C| &= |T| \quad	{\rm and}\cr
|\partial \widetilde C| &\le |\partial T|.}$$

\noindent {\bf Theorem 2.} If $C$ is a symmetric convex body in ${\Bbb
R}^n$ and $Q$ an $n$-dimensional cube then there is an affine image
$\widetilde C$ of $C$ satisfying

$$\eqalign{|\widetilde C| &= |Q|\quad {\rm and}\cr
|\partial \widetilde C| &\le |\partial Q|.}$$

Theorems 1 and 2 are proved in Section 1 below. The upper bounds for
surface area that are needed depend upon volume ratio estimates. The volume
ratio of an $n$-dimensional convex body $C$ is

$$vr(C) = \Big({|C|\over |{\cal E}|}\Big)^{1\over n}$$

\noindent where ${\cal E}$ is the ellipsoid of maximal volume included
in $C$. Similarly, if $X$ is an $n$-dimensional normed space, $vr(X)$ is
defined to be $vr(C)$ for any convex body $C$ which is the unit ball of $X$
in some representation of $X$ on ${\Bbb R}^n$. Section 2 of this paper
deals with a further question about volume ratios. It was proved in [B-M]
that the volume ratio of a finite-dimensional normed spaced can be bounded
above, solely in terms of the cotype-2 constant of the space: \ (see e.g.
[M-S] for definitions). In particular, finite-dimensional subspaces of
$L_1$ have uniformly bounded volume ratios. An isometric form of this
result is proved below: \ it is shown that for each $p\in [1,\infty],
\ell^n_p$ has maximal volume ratio among $n$-dimensional subspaces of
$L_p$.

This paper constitutes a sequel to the paper [B] which appeared recently.
\vfill\eject

\noindent {\bf \S 1. The reverse isoperimetric inequality.}

There are (at least) two ways in which to couch reverse forms of the
classical inequalities for convex bodies. In the case of the reverse
Santalo inequality, the expression to be estimated is $|C| |C^0|$. This
expression, the volume product for a symmetric convex body and its polar is
invariant under linear transformation of $C$. So there are bodies for which
the expression is minimal. Surface area does not behave well under linear
transformations. Although there are affine invariants which measure surface
area it seems natural to reverse the isoperimetric inequality by choosing
representatives of affine equivalence classes of bodies, as described in
the introduction: \ (for the affine invariant problem, see the appendix at
the end of this paper).

For many of the classical inequalities involving convex bodies, the
extremal bodies are ellipsoids or Euclidean balls and often this means that
the inequalities can be proved by well-known symmetrisation techniques. For
the reverse inequalities, one expects extremal bodies such as cubes or
tetrahedra: \ because of this, classical symmetrisation methods do not seem
to be readily applicable.

As was mentioned earlier, Theorems 1 and 2 are proved via volume ratio
estimates. Two well-known theorems of John [J], characterise ellipsoids of
maximal volume contained in convex bodies. These are stated here as lemmas.

\noindent {\bf Lemma 3.} Let $C$ be a symmetric convex body in ${\Bbb
R}^n$. The ellipsoid of maximal volume in $C$ is the Euclidean unit ball
$B^n_2$, if and only if $C$ contains $B^n_2$ and there is a sequence $(u_i)^m_1$
of contact points
between $B^n_2$ and $\partial C$ (i.e. unit vectors on the boundary of $C$)
and a sequence $(c_i)^m_1$ of positive numbers so that

$$\sum ^m_1 c_iu_i \otimes u_i = I_n.\eqno (1)$$

$\hfill \square$

Here, $u_i \otimes u_i$ is the rank-1 orthogonal projection onto the span
of $u_i$ and $I_n$ is the identity operator on ${\Bbb R}^n$. Condition (1)
shows that the $u_i$'s behave like an orthonormal basis to the extent that
for each $x\in {\Bbb R}^n$,

$$|x|^2 = \sum ^m_1 c_i\langle u_i, x\rangle^2.$$

\noindent The equality of the traces in (1) shows that

$$\sum ^m_1 c_i  = n.$$

\noindent {\bf Lemma 4.} Let $C$ be a convex body in ${\Bbb R}^n$ (not
necessarily symmetric). The ellipsoid of maximal volume in $C$ is $B^n_2$,
if and only if $C$ contains $B^n_2$ and
there are contact points $(u_i)^m_1$ and positive numbers
$(c_i)^m_1$ so that\medskip

\item{a)} $\sum\limits ^m_1 c_iu_i \otimes u_i = I_n$ \ \ and

\item{b)} $\sum\limits ^m_1 c_iu_i = 0$

$\hfill \square$

Theorems 1 and 2 are proved by combining the theorems of John with a
generalised convolution inequality of Brascamp and Lieb, [B-L]. A
``normalised'' form of this inequality was introduced in the author's
previous paper [B]. The normalisation is motivated by the theorems of John:
\ its principal advantage is that it ``automatically'' calculates the best
possible constant in the inequality. The theorem of Brascamp and Lieb is
stated here as a lemma.

\noindent {\bf Lemma 5.} Let $(u_i)^m_1$ be a sequence of unit vectors in
${\Bbb R}^n$ and $(c_i)^m_1$ a sequence of positive numbers so that

$$\sum^m_1 c_iu_i \otimes u_i = I_n.$$

\noindent For each $i$, let $f_i\colon \ {\Bbb R}\to [0,\infty)$ be
integrable. Then

$$\int_{{\Bbb R}^n} \prod ^m_{i=1} f_i(\langle
u_i,x\rangle)^{c_i} dx \le \prod ^m_{i=1} \bigg(\int_{\Bbb R} f
_i\bigg)^{c_i}.$$

$\hfill \square$

There is equality in Lemma 5 if the $f_i$'s are identical Gaussian
densities or if the $u_i$'s form an orthonormal basis of ${\Bbb R}^n$ (and
in some other cases). Lemma 5 is a generalisation of Young's convolution
inequality with best possible constant, proved independently by Beckner,
[Be].

It is relatively simple to combine Lemmas 3 and 5 to show that among
symmetric convex bodies in ${\Bbb R}^n$, the cube has exactly maximal
volume ratio. (This was done in [B] and the result is quoted below as the
$p=\infty$ case of Theorem 6.) Thus, if $C$ is a symmetric convex body
whose ellipsoid of maximal volume if $B^n_2$, then $|C| \le 2^n$. From this
it is easy to deduce Theorem 2.

\noindent {\bf Proof of Theorem 2.} Let $C$ be a symmetric convex body in
${\Bbb R}^n$. It is required to show that some affine image $\widetilde C$
of $C$ satisfies

$$|\partial \widetilde C| \le 2n |\widetilde  C|^{n-1\over n}$$

\noindent since these expressions are equal if $\widetilde C$ is a cube in
${\Bbb R}^n$. Choose $\widetilde C$ so that its ellipsoid of maximal volume
is $B^n_2$. By the remark above, $|\widetilde C| \le 2^n$. But, since
$B^n_2 \subset \widetilde C$,

$$\eqalignno{|\partial \widetilde C| &= \lim_{\varepsilon \to 0}
{|\widetilde C + \varepsilon B^n_2| - |\widetilde C|\over \varepsilon }
\le \lim_{\varepsilon \to 0} {|\widetilde C + \varepsilon \widetilde C| -
|\widetilde C|\over \varepsilon }
= \lim_{\varepsilon \to 0} |\widetilde C| \cdot {(1+\varepsilon)^n-1\over
\varepsilon }\cr
&= n|\widetilde C| = n|\widetilde C|^{n-1\over n} |\widetilde C|^{1\over
n}
\le 2n |\widetilde C|^{n-1\over n}.&\square}$$\medskip

In exactly the same way, Theorem 1 may be deduced from the volume ratio
estimate given by Theorem 1$'$ below. (In each case the argument loses
nothing, because each of the cube and tetrahedron has the property that all
its faces touch its ellipsoid if maximal volume.)

\noindent {\bf Theorem 1$'$.} Among all convex bodies in ${\Bbb R}^n$,
$n$-dimensional tetrahedra have maximal volume ratio.\medskip

\noindent {\bf Proof.} The $n$-dimensional regular tetrahedron that
circumscribes $B^n_2$ has volume

$${n^{n\over 2} (n+1)^{n+1\over 2}\over n!}.$$

\noindent So, it suffices to prove that if $C$ is a convex body whose
maximal ellipsoid is $B^n_2$ then

$$|C| \le {n^{n\over 2}(n+1)^{n+1\over 2}\over n!}.$$

\noindent By Lemma 4 there are unit vectors $(u_i)^m_1$ on $\partial C$ and
positive numbers $(c_i)^m_1$ so that

$$\eqalignno{&\sum ^m_1 c_iu_i \otimes u_i = I_n \quad {\rm and}&(2)\cr
&\sum^m_1 c_iu_i = 0.&(3)}$$

\noindent Since the $u_i$'s are contact points of $B^n_2$ and $\partial C$,

$$C\subset \{x\in {\Bbb R}^n\colon \ \langle x,u_i\rangle \le 1,  1
\le i \le m\} = K\qquad {\rm (say).}$$

\noindent It will be shown that $K$ has volume no more than

$${n^{n\over 2} (n+1)^{n+1\over 2}\over n!}.$$

Now, there is equality in Lemma 5 if the vectors  appearing in its
statement are orthogonal. The key to the following estimate is the
construction of a new sequence of vectors $(v_i)^m_1$ in ${\Bbb R}^{n+1}$
which would be orthogonal in the extreme case in which $K$ is a regular
tetrahedron. The estimate follows from an application of Lemma 5 to a
family of functions whose product is supported on a cone in ${\Bbb
R}^{n+1}$ whose cross-sections are similar to $K$.

Regard ${\Bbb R}^{n+1}$ as ${\Bbb R}^n\times {\Bbb R}$. For each $i$, let

$$\eqalignno{v_i &= \sqrt{n\over n+1} \Big(-u_i, {1\over \sqrt{n}}\Big) \in
{\Bbb R}^{n+1}\cr
\noalign{\hbox{and}}
d_i &= {n+1\over n} c_i.}$$

\noindent Then, for each $i, v_i$ is a unit vector and the identities (2)
and (3) above, together ensure that

$$\sum ^m_1 d_iv_i \otimes v_i = I_{n+1}.$$

\noindent Also, $\sum\limits ^m_1 d_i = n+1$. For each $i$, define
$f_i\colon \ {\Bbb R} \to [0,\infty)$ by

$$f_i(t) = \left\{\matrix{e^{-t}\hfill&{\rm if}&t\ge 0\hfill\cr
0\hfill&{\rm if}&t<0.\hfill}\right.$$

\noindent Finally, for $x\in {\Bbb R}^{n+1}$ set

$$F(x) = \prod ^m_{i=1} f_i(\langle v_i,x\rangle )^{d_i}.$$

\noindent By Lemma 5,

$$\int_{{\Bbb R}^{n+1}} F(x)dx \le \prod ^m_{i=1} \bigg( \int_{\Bbb R}
f_i\bigg)^{d_i}  = 1.\eqno (4)$$

Now, suppose $x = (y,r) \in {\Bbb R}^n \times {\Bbb R}$. For each $i$,

$$\langle v_i,x\rangle = {r\over \sqrt{n+1}} - \sqrt{n\over n+1} \langle
u_i, y\rangle.$$

\noindent Since $\sum\limits ^m_1 c_iu_i = 0$, there is some $j$ (depending
upon $y$) for which $\langle u_j,y\rangle \ge 0$. Hence, if $r<0, \langle
v_j,x\rangle < 0$ for some $j$ and so $F(x) = 0$. On the other hand, if
$r\ge 0$ then $F(x)$ is non-zero precisely if for every $i$

$$\langle u_i, y\rangle \le {r\over \sqrt{n}};$$

\noindent  in this case

$$\eqalign{F(x) &= \exp \bigg(-\sum ^m_1 d_i \Big[{r\over \sqrt{n+1}} -
\sqrt{n\over n+1} \langle u_i, y\rangle \Big]\bigg)\cr
&= \exp \bigg(-\sqrt{n+1}\cdot r + \sqrt{n+1\over n} \bigg\langle \sum^m_1
c_iu_i, y\bigg\rangle \bigg)\cr
&= \exp (-\sqrt{n+1} \cdot r).}$$

\noindent Thus, for each $r\ge 0$, the integral of $F$ over the hyperplane
$\{x\colon \ x_{n+1} = r\}$ is

$$e^{-\sqrt{n+1}r} \Big|{r\over \sqrt{n}} K\Big| = e^{-\sqrt{n+1}r}
\Big({r\over \sqrt{n}}\Big)^n |K|.$$

\noindent Therefore, from (4),

$$1 \ge |K| \int ^\infty_0 e^{-\sqrt{n+1}r} \Big({r\over \sqrt{n}}\Big)^n
dr = {|K|\cdot n!\over n^{n\over 2}(n+1)^{n+1\over 2}}$$

\noindent as required.$\hfill \square$\medskip

\noindent {\bf Remark.} It is possible to obtain estimates on outer volume
ratio, which involves the ellipsoid of minimal volume containing a body
$C$, by using a reverse form of the Brascamp and Lieb inequality. Among
symmetric convex bodies, the $n$-dimensional ``octahedron'' is extremal,
and among all convex bodies, the tetrahedron.\vfill\eject

\noindent {\bf \S 2. Volume ratios of subspaces of $L_p$.}

In [B], the author showed that among $n$-dimensional normed spaces,
$\ell^n_\infty$ has maximal volume ratio: \ (this fact was used in the
proof of Theorem 2 above). The inequality of Brascamp and Lieb can be
applied equally well to estimate the volume ratios of subspaces of $L_p$
for $1\le p < \infty$, since the volumes of their unit balls can be easily
expressed as convolutions: \ (this is done in Lemma 7 which was observed in
[M-P]).\medskip

\noindent {\bf Theorem 6.} Let $1 \le p \le \infty$. Among $n$-dimensional
subspaces of $L_p$, the space $\ell^n_p$ has exactly maximal volume ratio.
\medskip

The remainder of this section is devoted to the proof of Theorem 6.
\medskip

\noindent {\bf Lemma 7.} Let $K$ be a symmetric convex body in ${\Bbb R}^n$
with Minkowski gauge $\|\cdot\|$ and $1 \le  p < \infty$. Then

$$|K| = {1\over \Gamma(1+{n\over p})} \int _{{\Bbb R}^n} e^{-\|x\|^p} dx.$$

\noindent {\bf Proof.} Let $\sigma$ be the rotation invariant probability
on the Euclidean sphere $S^{n-1}$ and $v_n$ the volume of the Euclidean
unit ball $B^n_2$. Then

$$\eqalignno{\int _{{\Bbb R}^n} e^{-\|x\|^p} dx &= nv_n \int_{S^{n-1}}
\int ^\infty_0 e^{-\|r \theta\|^p} r^{n-1} dr d\sigma(\theta)\cr
&= nv_n \int_{S^{n-1}} \|\theta\|^{-n} d\sigma(\theta) \cdot \int ^\infty_0
e^{-r^p} r^{n-1} dr\cr
&= n|K| \int ^\infty_0 e^{-r^p} r^{n-1}dr\cr
&=\Gamma\Big( 1 + {n\over p}\Big)|K|.&\square}$$

\noindent {\bf Proposition 8.} Let $(u_i)^m_1$ be unit vectors in ${\Bbb
R}^n$ and $(c_i)^m_1$ positive numbers satisfying

$$\sum^m_1 c_iu_i \otimes u_i = I_n,$$

\noindent $(\alpha_i)^m_1$ positive numbers and $1\le p < \infty$. For
each $x\in {\Bbb R}^n$ set

$$\|x\| = \bigg( \sum ^m_1 \alpha_i|\langle
u_i,x\rangle|^p\bigg)^{1\over p}.$$

\noindent Then if $K$ is the unit ball of the space $({\Bbb R}^n, \|\cdot
\|)$,

$$|K| \le {2^n\Gamma(1+{1\over p})^n\over \Gamma(1+{n\over p})} \prod
^m_{i=1} \Big({c_i\over \alpha_i}\Big)^{c_i\over p}.$$

\noindent  {\bf Remark:} \ For $1\le p \le 2$. Meyer and Pajor, [M-P],
proved that the largest $n$-dimensional sections of the unit ball of
$\ell^m_p$ are those spanned by  $n$ standard unit vectors. This result is
a consequence of Proposition 8 applied with $\alpha_i = c_i^{p\over 2},
1 \le i \le m$.\medskip

\noindent {\bf Proof.} For $1\le i \le m$, let $f_i \colon \ {\Bbb R}\to
[0,\infty)$ be defined by

$$f_i(t) = \exp\Big(-{\alpha_i\over c_i} |t|^p\Big).$$

\noindent By Lemmas 7 and 5,

$$\eqalignno{|K| &= {1\over \Gamma(1+{n\over p})} \int _{{\Bbb R}^n}
\exp \bigg(-\sum \alpha _i |\langle u_i,x\rangle |^p\bigg)dx\cr
&= {1\over \Gamma(1+{n\over p})} \int _{{\Bbb R}^n} \prod^n_{i=1} \Big(\exp
{-\alpha_i\over c_i} |\langle u_i,x\rangle|^p\Big)^{c_i} dx\cr
&= {1\over \Gamma(1+{n\over p})} \int _{{\Bbb R}^n} \prod ^m_{i=1}
f_i(\langle u_i,x\rangle)^{c_i}dx\cr
&\le {1\over \Gamma(1+{n\over p})} \prod ^m_{i=1} \bigg(\int_{\Bbb R}
f_i\bigg)^{c_i}\cr
&= {1\over \Gamma(1+{n\over p})} \prod ^m_{i=1} \Big(2\Big({c_i\over
\alpha_i}\Big)^{1\over p} \Gamma\Big( 1 + {1\over
p}\Big)\Big)^{c_i}\cr
&= {2^n\Gamma(1+{1\over p})^n\over \Gamma(1+{n\over p})} \prod
^m_{i=1} \Big({c_i\over \alpha_i}\Big)^{c_i\over p}.&\square}$$

The deduction of Theorem 6 from the preceding proposition uses an important
lemma of Lewis [L], which extends John's theorem, Lemma 3.

\noindent {\bf Lemma 9.} Let $1\le p < \infty$ and $X$ be an
$n$-dimensional subspace of $\ell^m_p$. Then $X$ may be represented on
${\Bbb R}^n$ with norm given by

$$\|x\| = \bigg( \sum ^m_1 c_i|\langle u_i,x\rangle|^p\bigg)^{1\over
p}$$

\noindent where $(u_i)^m_1$ is a sequence of unit vectors and $(c_i)^m_1$ a
sequence of positive numbers satisfying

$$\sum^m_1 c_iu_i\otimes u_i = I_n.$$

$\hfill \square$

\noindent {\bf Remark.} The key point in Lemma 9 is that the same $c_i$'s
appear in both expressions. Thus, Proposition 8 may be applied with
$\alpha_i = c_i, 1 \le i\le m$; in this case the estimate of Proposition 8
has  a particularly simple form.

\noindent {\bf Proof of Theorem 6.} It may be assumed that $p<\infty$ and
it suffices to estimate the volume ratios of $n$-dimensional subspaces of
$\ell^m_p$ for each integer $m\ge n$. Let $X$ be such a space and assume
that it is represented on ${\Bbb R}^n$ in the way guaranteed by Lemma 9,
with unit ball $K$ (say). By Proposition 8,

$$|K| \le {2^n\Gamma(1+{1\over p})^n\over \Gamma(1+{n\over p})}.$$

\noindent The latter expression is the volume of the unit ball of
$\ell^n_p$ in its normal representation on ${\Bbb R}^n$. So it is enough to
check that $K$ contains an Euclidean ball of radius

$$\left\{\matrix{n^{{1\over 2}-{1\over p}}\hfill&{\rm if}&1 \le p \le
2\hfill\cr
1\hfill&{\rm if}&p>2\hfill.}\right.$$

\noindent In the first case, for every $x\in {\Bbb R}^n$,

$$\eqalign{\|x\|^p &= \sum ^m_1 c_i |\langle u_i,x\rangle|^p\cr
&\le \bigg(\sum ^m_1 c_i\bigg)^{1-{p\over 2}} \bigg(\sum ^m_1 c_i \langle
u_i,x\rangle^2\bigg)^{p\over 2}\cr
&= n^{1-{p\over 2}} |x|^p}$$

\noindent ($|x|$ being the Euclidean length of $x$). In the second case,

$$\|x\|^p = \sum ^m_1 c_i |\langle u_i,x\rangle|^p
\le \sum ^m_1 c_i |x|^{p-2} \langle u_i,x\rangle^2
= |x|^{p-2}\cdot |x|^2 = |x|^p.$$

$\hfill \square$\vfill\eject

\noindent {\bf Appendix.} For a convex body $C$ in ${\Bbb R}^n$ and unit
vector $\theta \in {\Bbb R}^n$, let $P_\theta C$ be the orthogonal
projection of $C$ onto the 1-codimensional subspace of ${\Bbb R}^n$
perpendicular to $\theta$. It was observed by Petty in [P$_1$], that the
expression

$$\bigg( |C|^{n-1} \int_{S^{n-1}} |P_\theta C|^{-n} d\sigma(\theta
)\bigg)^{-{1\over n}}\eqno (5)$$

\noindent is invariant under invertible affine transformations of the body
$C$. This expression measures ``minimal surface area'' in the sense of
Theorems 1 and 2. To see this, note that, on the one hand,
the Cauchy formula for surface area states that for a
convex body $C$ in ${\Bbb R}^n$,

$$|\partial C| = {nv_n\over v_{n-1}} \int_{S^{n-1}} |P_\theta
C|d\sigma(\theta)$$

\noindent and so by H\"older's inequality

$${|\partial C|\over |C|^{n-1\over n}} \ge {nv_n\over v_{n-1}}
\bigg(|C|^{n-1} \int_{S^{n-1}} |P_\theta C|^{-n} d\sigma
(\theta)\bigg)^{-{1\over n}}.$$

\noindent On the other hand, Theorem 6 (for $p=1$) states that if $X$ is an
$n$-dimensional subspace of $L_1$ then

$$vr(X) \le \Big({2^n\Gamma(1+{n\over 2})\over \Gamma(1+n)\pi^{n\over
2}}\Big)^{1\over n}\le \sqrt{2e\over \pi}.$$

\noindent From this it is easy to deduce that each convex body $C$ has an
affine image  $\widetilde C$ for which

$$\eqalign{{|\partial \widetilde C|\over |\widetilde C|^{n-1\over n}} &\le
\sqrt{2e\over \pi} {nv_n\over v_{n-1}} \bigg(|\widetilde C|^{n-1}
\int_{S^{n-1}} |P_\theta \widetilde  C|^{-n}
d\sigma(\theta)\bigg)^{-{1\over n}}\cr
&= \sqrt{2e\over \pi} {nv_n\over v_{n-1}} \bigg(|C|^{n-1} \int _{S^{n-1}}
|P_\theta C|^{-n} d\sigma(\theta))^{-{1\over n}}.}$$

A strong isoperimetric inequality of Petty, [P$_2$], states that  the
expression (5) is minimised by the Euclidean balls. It would be possible
to reverse the isoperimetric inequality by determining the  bodies which
maximise (5). It seems likely that the cube and tetrahedron are extremal
for this modified problem: \ this is certainly true if $n=2$.
\vfill\eject

\noindent {\bf References}

\item{[B]} K.M. Ball, Volumes of sections of cubes and related problems,
Israel seminar on Geometric Aspects of Functional Analysis (J.
Lindenstrauss
and V.D. Milman eds.), Lecture Notes in Mathematics \#1376,
Springer-Verlag, (1989), 251-260.

\item{[Be]} W. Beckner, Inequalities in Fourier analysis, Annals of Math.
102 (1975), 159-182.

\item{[B-M]} J. Bourgain and V.D. Milman, New volume ratio properties for
convex symmetric bodies in ${\Bbb R}^n$, Inventiones Math. 88 (1987),
319-340.

\item{[B-L]} H.J. Brascamp and E.H. Lieb, Best contants in Young's
Inequality, its converse, and its generalization to more than three
functions, Advances in Math. 20 (1976), 151-173.

\item{[J]} F. John, Extremum problems with inequalities as subsidiary
conditions, Courant Anniversary Volume, Interscience, New York (1948),
187-204.

\item{[L]} D.R. Lewis, Finite dimensional subspaces of $L_p$, Studia Math.
63 (1978), 207-212.

\item{[M-P]} M. Meyer and A. Pajor, Sections of the unit ball of
$\ell^n_p$, Journal of Funct. Anal. 80 (1988), 109-123.

\item{[M]} V.D. Milman, An inverse form of the Brunn-Minkowski inequality
with applications to local theory of normed spaces, C.R. Acad. Sci. Paris,
302 (Ser. 1), \#1 (1987).

\item{[M-S]} V.D. Milman and G. Schechtman, Asymptotic Theory of Finite
Dimensional Normed Spaces, Lecture notes in mathematics \#1200,
Springer-Verlag, Berlin-Heidelberg-New York, (1985).

\item{[P$_1$]} C.M. Petty, Projection bodies, Proc. Colloq. on Convexity,
Copenhagen, 1967, 234-241.

\item{[P$_2$]} C.M. Petty, Isoperimetric problems, Proc. Conf. on Convexity
and Combinatorial Geometry (Univ. of Oklahoma, June 1971), 1972, 26-41.

\end